\documentclass[preprint,sort&compress,12pt]{elsarticle}
\usepackage[pagewise]{lineno}
\usepackage{bbm}
\usepackage{amstext}
\usepackage{amsmath}
\usepackage{amssymb}
\usepackage{mathrsfs}
\usepackage{stmaryrd}
\usepackage{bm}
\usepackage{pstricks}
\usepackage{pst-coil}
\usepackage{pst-3d}
\usepackage{amsfonts}
\usepackage{amsthm}
\usepackage{CJK}
\allowdisplaybreaks

\newcommand{\mf}{\mathbf}
\newcommand{\mm}{\mathrm}

\newcommand{\be}{\begin{equation}}
\newcommand{\bea}{\begin{equation}\begin{aligned}}
\newcommand{\beas}{\begin{equation*}\begin{aligned}}
\newcommand{\eeas}{\end{aligned}\end{equation*}}
\newcommand{\eea}{\end{aligned}\end{equation}}
\newcommand{\ee}{\end{equation}}

\renewcommand{\div}{{\rm div }}

\begin{document}
\begin{CJK*}{GBK}{song}
	\begin{frontmatter}
		\title{Small Alfv\'en Number Limit for the Global-in-time Solutions of Incompressible MHD Equations with General
Initial Data}
\author[Yuan]{Yuan Cai}
		\ead{caiy@fudan.edu.cn}
\author[Cui]{Xiufang Cui}
		\ead{cuixf@lzu.edu.cn} 
		\author[mmL1]{Fei Jiang}
		\ead{jiangfei0591@163.com}
		\author[mmL1]{Hao Liu}
		\ead{liuh20230111@163.com}
\address[Yuan]{School of Mathematical Sciences, Fudan University, Shanghai 200433, China.}
\address[Cui]{School of Mathematics and Statistics, Lanzhou University,  Lanzhou, 730000,  China.}   
		\address[mmL1]{School of Mathematics and Statistics, Fuzhou University, Fuzhou, 350108, China.}				

\begin{abstract}
The small Alfv\'en number (denoted by $\varepsilon$) limit (one type of large parameter limits, i.e. singular limits) in magnetohydrodynamic (abbr. MHD) equations was first proposed by Klainerman--Majda in (Comm. Pure Appl. Math. 34: 481--524, 1981). Recently Ju--Wang--Xu mathematically verified that the \emph{local-in-time} solutions of three-dimensional (abbr. 3D)  ideal (i.e. the absence of the dissipative terms) incompressible MHD equations with general initial data in $\mathbb{T}^3$ (i.e. a spatially periodic domain) tend to a solution of 2D ideal MHD equations in the distribution sense as $\varepsilon\to 0$ by Schochet's fast averaging method in (J. Differential Equations, 114: 476--512, 1994). In this paper, we revisit the small Alfv\'en number limit in $\mathbb{R}^n$ with $n=2$, $3$, and develop another approach, motivated by Cai--Lei's energy method in (Arch. Ration. Mech. Anal. 228: 969--993, 2018), to establish a new conclusion that the \emph{global-in-time} solutions of incompressible MHD equations (including the viscous resistive case) with general initial data converge to zero  as $\varepsilon\to 0$ for any given time-space variable $(x,t)$ with $t>0$. 
In addition, we find that the large perturbation solutions and vanishing phenomenon of the nonlinear interactions also exist in the \emph{viscous resistive} MHD equations for small Alfv\'en numbers, and thus extend Bardos et al.'s  results of the \emph{ideal} MHD equations in (Trans Am Math Soc 305: 175--191, 1988). 
\end{abstract}
\end{frontmatter}
\newtheorem{thm}{Theorem}[section]
\newtheorem{lem}{Lemma}[section]
\newtheorem{pro}{Proposition}[section]
\newtheorem{concl}{Conclusion}[section]
\newtheorem{cor}{Corollary}[section]
\newproof{pf}{Proof}
\newdefinition{rem}{Remark}[section]
\newtheorem{definition}{Definition}[section]
	
\section{Introduction}\label{sect1}
\numberwithin{equation}{section}
In this paper, we investigate the small Alfv\'en number limit for the global-in-time solutions of the following incompressible MHD system in $\mathbb{R}^n$ with $n=2$, $3$:
	\begin{equation}\label{1.1}
	\begin{cases}
	 \rho v_t+  \rho v\cdot\nabla v- \mu\Delta v+\nabla P =\lambda  M\cdot\nabla M/4\pi ,\\[1mm]
	M_t+v\cdot\nabla M-\nu\Delta M=M\cdot\nabla v,\\[1mm]
	\mathrm{div}v=\mathrm{div}{M}=0.
	\end{cases}
	\end{equation}
	The unknown functions ${v}:= {v}(x,t)$, $M:= {M}(x,t)$ and $P:= P(x,t)$ denote the velocity, the magnetic field and the sum of kinetic  pressure and magnetic pressure, resp.. The nonnegative constants $\rho>0$, $\mu\geqslant $ and $\nu\geqslant 0$ represent the  coefficients of density, viscosity and resistivity, resp.. In addition, $\lambda$ is the magnetic permeability. The well-posedness problem of the above system have  been widely investigated, please refer to \cite{MR3780143,MR4641656,MR3738384,MR716200,MR3678491,MR3666563,MR3590662} and the references cited therein. 

\subsection{Motivation}
Now we consider a rest state $(0,\bar{M})$ for the system \eqref{1.1}, and then define
$$H:=\sqrt{\frac{\lambda}{4\pi\rho}}M-\bar{M},\ p:=\frac{P}{\rho},\  \bar{M}:=\sqrt{ \frac{\lambda}{4\pi\rho}}m\mf{e}^{1},\ \tilde{\mu}:=\frac{\mu}{\rho}\mbox{ and }\varepsilon:=\sqrt{\frac{4\pi\rho}{\lambda m^2}},$$
 where $\bar{M} $ is an impressed magnetic field,  $\mathbf{e}^1\in \mathbb{R}^n$ represents the unit vector with the first component being $1$, $\lambda$ denotes the permeability of vacuum and $\varepsilon$ is called the \emph{Alfv\'en number} in this paper. Thus we have the following perturbation system for $(v,H)$:
\begin{equation}\label{1sa.1}
\begin{cases}
v_t+  v\cdot\nabla v- \tilde{\mu}\Delta v+\nabla p =(\varepsilon^{-1} \mathbf{e}^1+ H )\cdot \nabla H, \\[1mm]
H_t+v\cdot\nabla H-\nu\Delta H=(\varepsilon^{-1}  \mf{e}^{1} +H)\cdot\nabla v,\\[1mm]
\mathrm{div}v=\mathrm{div}H=0.
\end{cases}
\end{equation}
Next we force on the relevant progress for the small Alfv\'en number limit problem.
 
Small Alfv\'en number limit (one type of large parameter limits  \cite{MR871107,MR665380,MR957005}, i.e. singular limits, due to $m\to \infty$ as $\varepsilon\to 0$) is one of the distinguished singular limits for MHD \cite{MR3803773}. For 3D compressible ideal MHD system, the singular problem with respect to the Alfv\'en number in $\mathbb{T}^3$ (i.e. a spatially periodic domain) was first proposed by Klainerman--Majda \cite{MR615627}, in which the small Mach number limit was also investigated. 
However, they cannot prove convergence to an appropriate reduced system. To obtain the limit system,
Browning--Kreiss posed more unnatural assumptions on initial data, i.e., high-order time derivatives
of the solution are assumed to be bounded uniformly with respect to the small Alfv\'en numbers at $t=0$ \cite{MR665380}.

Later Goto forced on the incompressible case of the 3D ideal MHD system (without viscosity) under the natural assumptions on the initial data and first determined that the limiting system becomes essentially the system of the 2D motion \cite{MR1039472}.  Rubino further considered the viscous case \cite{MR1339828}. Recently Jiang--Ju--Xu extended Goto's result to the case that the domain is a slab with horizontal periodicity, and found that the limiting system is a 2D Euler system coupled with a linear transport equation \cite{MR4028271}. 
Similar results can be found in the corresponding compressible case, see \cite{MR3836806,MR3962511,MR4221327,MR4853427,MR4425021} for examples.  
 
It should be emphasized that the above works with a limit system  are only valid for well-prepared initial data. When the initial data is ill-prepared, the fast oscillating waves will be developed and the mathematical analysis is more complicated. 
Recently Ju--Wang--Xu mathematically verified that the \emph{local-in-time} solutions of 3D incompressible  ideal MHD equations with general initial data in $\mathbb{T}^3$ tend to a solution of 2D ideal MHD equations in the distribution sense as $\varepsilon\to 0$ \cite{MR4645635} by Schochet's fast averaging method in \cite{MR1303036}.
However, the small Alfv\'en number limit in 3D \emph{compressible ideal} MHD equations with \emph{general initial data}  is one of the largely wide open problems (see Majda's remark on pp. 72 in his monograph \cite{MR748308}). 
We mention the small Alfv\'en number limit of \emph{global-in-time} weak solutions of 3D non-isentropic compressible \emph{viscous resistive}  MHD equations has been investigated by Kuku$\mathrm{\check{c}}$ka \cite{MR2805860}, however the limit system is still 3D. 

In this paper, we revisit the small Alfv\'en number limit in $\mathbb{R}^n$ with $n=2$, $3$ for the incompressible MHD equations  with  general initial data for the case
\begin{equation}
\label{2504011059}
\tilde{\mu}=\nu \geqslant 0,
\end{equation} and develop another approach, motivated by Cai--Lei's energy method in \cite{MR3780143}, to establish a new conclusion that the \emph{global-in-time} classical solutions   converge to zero  as $\varepsilon\to 0$ for any given time-space variable $(x,t)$ with $t>0$, please refer to the decay of solutions with respect to $\varepsilon$ in Theorem \ref{thsm3} for details. 

\subsection{Reformulation}	
To conveniently investigate the small Alfv\'en number limit, we shall first rewrite the system \eqref{1sa.1}. Let us recall the Els$\ddot{\rm a}$sser variables  \cite{MR920153} 
\begin{align}
\label{2503222225}
{\Lambda}_+:=v+H\mbox{ and }{\Lambda}_-:=v-H,
\end{align}
which allow us to rewrite the perturbation equations with \eqref{2504011059} for $(\Lambda_+,\Lambda_-, p)$ as follows:
\begin{equation}\label{1.3}
\begin{cases}
\partial_t{\Lambda}_+-\varepsilon^{-1}  \partial_1 {\Lambda}_++ {\Lambda}_- \cdot\nabla {\Lambda}_++\nabla p= \nu\Delta {\Lambda}_+ ,\\[1mm]
\partial_t\Lambda_-+\varepsilon^{-1} \partial_1 {\Lambda}_-+ {\Lambda}_+ \cdot\nabla {\Lambda}_-+\nabla p= \nu\Delta {\Lambda}_-,\\[1mm]
\mathrm{div}\Lambda_+=\mathrm{div} \Lambda_-=0.
\end{cases}
\end{equation}

We further use the variables
\begin{align*}
t_*:=\varepsilon ^{-1} t,\  \Lambda_*^{\pm}:={\Lambda_\pm}(x,\varepsilon{t_*})\mbox{ and }p_*:= \varepsilon p(x,\varepsilon{t_*} ) 
\end{align*}
to rewrite the system \eqref{1.3} for $\Lambda_*^\pm$ and $p_*$, and then omit the subscript $*$ to obtain the following Cauchy problem 
\begin{equation}\label{1.4}
\begin{cases}
\Lambda_t^\pm\mp\partial_1 {\Lambda}^\pm+\varepsilon  {\Lambda}^\mp \cdot\nabla {\Lambda}^\pm+\nabla p= \varepsilon\nu\Delta {\Lambda}^\pm ,\\[1mm]
\mathrm{div}{\Lambda}^\pm=0,\\
{\Lambda}^\pm|_{t=0}={\Lambda}_\pm^0.
\end{cases}
\end{equation} 
We mention that we can also directly establish our main result in Theorem \ref{thsm3} from the system \eqref{1.3}, but the derivations are extremely complicated. 

In next section, we will present the existence of global-in-time solutions $\Lambda^\pm$ of the Cauchy problem \eqref{1.4} with a uniform energy estimate with  respect to Alfv\'en {number} in Theorem \ref{thm1}, the error estimate between the (nonlinear) solution of  the Cauchy problem  \eqref{1.4} and the solution 
 of the corresponding linear system in Theorem \ref{thm2}, 
and the decay behavior of solutions $\Lambda^\pm$ in Theorem \ref{thm3} in sequence. We mention that the error estimate in Theorems \ref{thm2} provides the decay-in-$\varepsilon$ of solutions for the case $\nu=0$ in Theorem \ref{thm3}. Then we easily obtain the Alfv\'en number limit behavior of $(v,H)$ in Theorem \ref{thsm3} from  Theorem \ref{thm3}. In addition, since the energy estimate of solutions in Theorem \ref{thm1} is also uniformly bound with  respect to dissipative coefficient $\nu$,  we also provide an additional result for the non-dissipative limit in Theorem \ref{2504160948}. 

\section{Main results}
Before stating the main results, we first introduce some notations used frequently throughout this paper.

(1) Basic notations: $\mathbb{R}_{0}^{+}:=[0,\infty)$, $\langle  \sigma\rangle:= \sqrt{1+|\sigma|^2}$, $\int:=\int_{\mathbb{R}^{n}}$, $\alpha:=\left(\alpha_{1},  \cdots,\alpha_{n}\right)$ is the multi-index with respect to the variable $x$, i.e. $\partial^\alpha=\partial_{x_1}^{\alpha_1}\cdots\partial_{x_n}^{\alpha_n}$.  $c$  represents a generic positive constant, which at most depends on the parameters $n$, $k$ and ${s}$ appearing in Theorem \ref{thm1},  and may vary from line to line.  
	 $A \lesssim B$ means that  $A \leqslant c B$. A pseudo-differential operator $|\nabla|^{-1}$ is defined by $|\nabla|^{-1}f:=\mathcal{F}^{-1}_{\xi\to x}\left(|\xi|^{-1}\mathcal{F}_{x\to\xi}(f(x))\right)$, where $\mathcal{F}$ and $\mathcal{F}^{-1}$ are the Fourier transform and inverse Fourier transform, resp.. In addition,
\begin{align*}
\mm{sign}(\nu):=
\begin{cases}
0&\mbox{for }\nu=0;\\
1&\mbox{for }\nu> 0.
\end{cases}
\end{align*}  
	
	(2) Simplified norms and function spaces:
	\begin{align*}
	&L^2:=L^2(\mathbb{R}^n),\ H^{i}:=W^{i, 2}\left(\mathbb{R}^{n}\right), \ H^i_{\sigma}:=\{w\in H^i~|~\mm{div}w=0\},\\
&\|\cdot\|_i:=\|\cdot\|_{H^i},\ \|(f,g)\|_{i}:=\sqrt{\|f\|_i^2+\|g\|_i^2},\ \|\cdot\|_{L^p_tL^q_x}:=
\|\cdot\|_{L^p((0,t),L^q)},
\end{align*}
where  $i \geqslant 0$ is integer, and $1 \leqslant p$, $q\leqslant \infty$.

 (3) Energy/dissipation functionals:
\begin{align*}
{\mathcal{E}}_{k}(t):=&\mm{sign}(\nu)\left\|\left({|\nabla|^{-1}}
\Lambda^+,{|\nabla|^{-1}}
\Lambda^-\right)\right\|_0^2+\left\|\left(\langle x+ \mathbf{e}^1t\rangle^s\Lambda^+,\langle x- \mathbf{e}^1t\rangle^s\Lambda^-\right)\right\|_0^2\\ 
& +\sum_{|\alpha|=1}^{ k}\left\|\left(\langle x+ \mathbf{e}^1t\rangle^{2s}\partial^\alpha\Lambda^+,\langle x- \mathbf{e}^1t\rangle^{2s}\partial^\alpha\Lambda^- \right)\right\|_0^2,\\ 
\mathcal{W}_{k}(t):=&\left\|\left(\frac{\langle x+ \mathbf{e}^1t\rangle^s\Lambda^+}{\langle x_1- t\rangle^{s}},\frac{\langle x- \mathbf{e}^1t\rangle^s\Lambda^-}{\langle x_1+ t\rangle^{s}}\right)\right\|_0^2 \\&+\sum_{|\alpha|=1}^{k}\left\|\left(\frac{\langle x+ \mathbf{e}^1t\rangle^{2s}\partial^\alpha\Lambda^+}{\langle x_1- t\rangle^{s}},\frac{\langle x- \mathbf{e}^1t\rangle^{2s}\partial^\alpha\Lambda^-}{\langle x_1+ t\rangle^{s}}\right)\right\|_0^2\\
\mathcal{D}_{k}(t):=
	&\left\|\left(\Lambda^+,\Lambda^-\right)\right\|_0^2+
\left\|\left(\langle x+ \mathbf{e}^1t\rangle^s \nabla\Lambda^+,\langle x- \mathbf{e}^1t\rangle^s \nabla\Lambda^-\right)\right\|_0^2\\
 & +\sum_{|\alpha|=1}^{k}\left\|\left(\langle x+ \mathbf{e}^1t\rangle^{2s}\nabla \partial^\alpha\Lambda^+,\langle x- \mathbf{e}^1t\rangle^{2s}\nabla \partial^\alpha\Lambda^-\right)\right\|_0^2 .
	\end{align*} 
In addition, 
\begin{align*}
{\mathcal{E}}_{k}^0:=&
\mm{sign}(\nu)\left\|\left(|\nabla|^{-1}\Lambda^0_+,|\nabla|^{-1}\Lambda^0_-\right)\right\|_0^2+\left\|\langle x\rangle^s\left(\Lambda^0_+, \Lambda^0_-\right)\right\|_0^2\\
&+\sum_{|\alpha|=1}^{k} 
\left\|\langle x\rangle^{2s}\left(\partial^\alpha\Lambda^0_+, \partial^\alpha\Lambda^0_-\right)\right\|_0^2. \end{align*}

\subsection{Results for the transformed system \eqref{1.4}}
We introduce the first result concerning the existence of global-in-time solutions to the Cauchy problem \eqref{1.4}  with  a uniform estimate with respect to  the Alfv\'en number $\varepsilon $.
\begin{thm}\label{thm1}
Let $n=2$, $3$, $\nu\geqslant 0$,  $\varepsilon\nu\leqslant 1/2$, $1<2s <4/{3}$,
$\Lambda_\pm^{0} \in H^k_\sigma$  with the integer $k\geqslant4$ and 
\begin{align}
\label{2022504161839}
\gamma\leqslant 2\mbox{ with }\varepsilon\in (0,1] \mbox{ or }\gamma= 2\mbox{ with }\varepsilon>0.
\end{align}
There {exist a small constant}  $c_1\in(0,1)$ and large constant $c_2$, which depend only on $n$, $k$ and ${s}$ such that, for any $\varepsilon$ satisfying
\begin{equation}\label{1.6}
{\mathcal{E}}_{k}^0  \leqslant  c_1 \varepsilon^{-\gamma},
\end{equation}
 the Cauchy problem \eqref{1.4} defined on $\mathbb{R}^n\times \mathbb{R}_0^+$ admits a unique global-in-time classical solution $\Lambda^\pm $, which satisfy
 		\begin{align}
 \label{2504160945}
  \mathcal{E}_k (t)+\int_0^t(\mathcal{W}_k(\tau)+\varepsilon\nu \mathcal{D}_k(\tau))\mm{d}\tau\leqslant c_2  \mathcal{E}^0_k\mbox{ for a.e. }t>0.
 		\end{align}
\end{thm}
\begin{rem}  For $\gamma\in (0,2]$, we call \eqref{1.6} the condition of small Alfv\'en number.
\end{rem}

Now we state the second result concerning the error estimate between the (nonlinear) solution of  the Cauchy problem  \eqref{1.4} and the solution 
 of the corresponding linear problem in Theorem \ref{thm2}. 
\begin{thm}\label{thm2}
Let $\Lambda^\pm$ be the global solutions in Theorem \ref{thm1}, $\Lambda^\pm_{\mm{N},\mm{L}}=\Lambda^\pm-\Lambda^\pm_\mm{L}$ and 
$\Lambda^\pm_\mm{L}$ the unique classical solutions to the linear (pressureless) problem in $\mathbb{R}^n\times \mathbb{R}_0^+$:
 \begin{align}\label{1.11}
 &\begin{cases}
 \partial_t\Lambda^\pm_\mm{L}\mp\partial_1\Lambda^\pm_\mm{L}-\varepsilon\nu\Delta {\Lambda}_\mm{L}^\pm=0,\\[1mm]  
 \div \Lambda^\pm_\mm{L}=0,\\[1mm]
\Lambda^\pm_\mm{L}|_{t=0}=\Lambda_\pm^0=\Lambda^\pm|_{t=0}.
 \end{cases} 
 \end{align} 
If the initial data additionally satisfy
\begin{align} 
\label{2504221950}
& \sum_{|\alpha|=k} 
\left\|\langle x\rangle^{3s}\left(\partial^\alpha\Lambda^0_+, \partial^\alpha\Lambda^0_-\right)\right\|_0<\infty,  \end{align} 
 it holds that \begin{align}
&\mathcal{E}_{k-1}^{\mm{N},\mm{L}}(t)+\int_0^t({\mathcal{W}}_{k-1}^{\mm{N},\mm{L}}(\tau)+
\varepsilon\nu  {\mathcal{D}}_{k-1}^{\mm{N},\mm{L}}(\tau)){\mm{d}}\tau\nonumber \\
&\lesssim \varepsilon \mathcal{E}^0_k \left(\sqrt{\mathcal{E}^0_k} 
+\sum_{|\alpha|=k} 
\left\|\langle x\rangle^{3s}\left(\partial^\alpha\Lambda^0_+, \partial^\alpha\Lambda^0_-\right)\right\|_0 \right)\mbox{ for any }t>0,\label{1.12}
\end{align} 
where   $\mathcal{E}_{k-1}^{\mm{N},\mm{L}}(t)$, $\mathcal{W}_{k-1}^{\mm{N},\mm{L}}(t)$  resp. $\mathcal{D}_{k-1}^{\mm{N},\mm{L}}(t)$ are defined by $\mathcal{E}_{k-1}(t)$, $\mathcal{W}_{k-1}(t)$ resp. $\mathcal{D}_{k-1}(t)$ with $\Lambda^\pm_{\mm{N},\mm{L}}$ in the place of $\Lambda^\pm$.
\end{thm}
 
Finally we present the decay behavior  of solutions in Theorem \ref{thm1}.
\begin{thm}\label{thm3}
Let $\Lambda^\pm$ be the global-in-time solutions in Theorem \ref{thm1}, it holds that
\begin{align}
|\Lambda^\pm(x,t)|
\lesssim 
 |x\pm \mathbf{e}^1t|^{-s}\sqrt{\mathcal{E}_k^0}
\mbox{ for any }(x,t)\in \mathbb{R}^n\times \mathbb{R}^+_0.
\label{2504140938}
\end{align}
If additionally the initial data have the compact supported in $B_R(0):=\left\{x\in\mathbb{R}^n\big{|}|x|<R\right\} $, the following estimate holds for the case $\nu=0$
 \begin{align}\label{2504140939}
|\Lambda^\pm(x,t)|\lesssim \sqrt{\varepsilon} (\mathcal{E}_k^0)^{\frac{3}{4}}
\end{align}
for any $x\in B_R(0)$ and for any $t>2R$.
    \end{thm} 

Now we mention the proofs of Theorems \ref{thm1}--\ref{thm3}.	The existence of the global-in-time solutions with small data has been proved by  the first author and Z. Lei for the Cauchy problem \eqref{1.4} for $\varepsilon=1$  \cite{MR3780143} by using the inherent strong null structure of the systems introduced
by Lei \cite{MR3552009} and the ghost weight technique introduced by Alinhac \cite{MR1856402}. 
Following the first author and Lei's energy method in \cite{MR3552009},  we naturally take the following \emph{a priori}  assumption for solutions defined on $\mathbb{R}^n\times [0,T)$
 \begin{align}
  \sup_{0\leqslant t<T}\mathcal{E}_k (t)+\int_0^T(\mathcal{W}_k(\tau)+\varepsilon\nu \mathcal{D}_k(\tau))\mm{d}\tau\leqslant \iota {c}_1  \varepsilon^{-\gamma}\mbox{ for some constant }\iota>0,
  \label{2504171021}
 \end{align}
and then derive the following \emph{a priori} estimate of solutions 
 \begin{align}
  &  \sup_{0\leqslant t<T}\mathcal{E}_k ( t)+ \int_0^T\left((1- \varepsilon\nu)  \mathcal{W}_k ( \tau)+   \varepsilon\nu \mathcal{D}_k (  \tau) \right)\mm{d}\tau \nonumber \\
  &\leqslant  \tilde{c}_2 \left(\mathcal{E}_k^0 + \varepsilon  \sup_{0\leqslant t<T}\sqrt{\mathcal{E}_k(  t)  } \int_0^T \mathcal{W}_k (\tau) \mm{d}\tau\right) \leqslant \tilde{c}_2 \left(\mathcal{E}_k^0 + \varepsilon^{1-\frac{\gamma}{2}}\sqrt{ \iota {c}_1   } \int_0^T \mathcal{W}_k ( t) \mm{d}\tau\right),
  \label{2504162226}
 \end{align}
 where $\tilde{c}_2$ depends only on $n$, $k$ and ${s}$.
 Recalling $\nu \varepsilon\leqslant 1/2$, both cases in \eqref{2022504161839} and \eqref{1.6},  then choosing  $12 \tilde{c}_2\leqslant  \iota$ and sufficiently small  $ {c}_1$,  and finally we  infer from the above estimate that  \begin{align}
    \sup_{0\leqslant t<T}\mathcal{E}_k (   t)+\int_0^T(
   \mathcal{W}_k ( \tau)+  \varepsilon\nu  \mathcal{D}_k (  \tau) )\mm{d}\tau \leqslant 
6 \tilde{c}_2 \mathcal{E}_k^0\leqslant\iota{c}_1  \varepsilon^{-\gamma}/2. 
\label{2504161849}
 \end{align} 
The above \emph{a priori} stability estimate, together with a local-in-time well-posedness result on the Cauchy problem \eqref{1.4}, immediately yields Theorem \ref{thm1}.	We will provide the derivation of \eqref{2504162226} in Section \ref{sec2}.

We also follow the derivation of \eqref{2504161849} to establish Theorem \ref{thm2} in Section \ref{sec3}. 
Finally, exploiting the Sobolev embedding inequality $H^2\hookrightarrow L^\infty$ and the solution of the transport equation, we can easily deduce Theorem {\ref{thm3}} from Theorems \ref{thm1} and \ref{thm2} in Section \ref{sec4}.

\subsection{Results for the original system \eqref{1sa.1} with $\tilde{\mu}=\nu$}

Recalling the relations in \eqref{2503222225}, it is easy to see that 
$$v=({\Lambda}^+(x,\varepsilon ^{-1} t)+{\Lambda}^-(x,\varepsilon ^{-1} t))/2\mbox{ and }H:=({\Lambda}^+(x,\varepsilon ^{-1} t)-{\Lambda}^-(x,\varepsilon ^{-1} t))/2.$$
Consequently we immediately write the following conclusion, which includes the small Alfv\'en number limit, for the solution $(v,H)$ of the original system \eqref{1sa.1} from Theorems \ref{thm1}--\ref{thm3}.  
\begin{thm}\label{thsm3}
Let $n=2$, $3$, $\nu\geqslant 0$,  $\varepsilon\nu\leqslant 1/2$, $1<2s <4/{3}$,
$(v^{0}, H^0) \in H^k_\sigma$  with the integer $k\geqslant4$ and 
\begin{align*}
\gamma\leqslant 2\mbox{ with }\varepsilon\in (0,1] \mbox{ or }\gamma= 2\mbox{ with }\varepsilon>0.
\end{align*} 
There {exists a small constant}  $c_3\in(0,1)$, which depends only on $n$, $k$ and ${s}$ such that, for any $\varepsilon$ satisfying
\begin{equation*} 
{ {E}}_{}^0:=\mm{sign}(\nu)\left\||\nabla|^{-1}(v^{0},H^0) \right\|_0^2 +\left\|\langle x\rangle^s(v^{0},H^0) \right\|_0^2+\sum_{|\alpha|=1}^k\left\|\langle x\rangle^{2s}\partial^\alpha(v^{0},H^0) \right\|_0^2\leqslant  c_3 \varepsilon^{-\gamma},
\end{equation*}
 the Cauchy problem of the system \eqref{1sa.1} defined on $\mathbb{R}^n\times \mathbb{R}_0^+$ admits a unique global-in-time classical solution $(v,H)$, which enjoys the following properties:
\begin{enumerate}[(1)]
  \item {the vanishing phenomenon of the nonlinear interactions under the additional regularity of initial data 
 \begin{align}
 \label{25042219501}
& \sum_{|\alpha|=k} 
\left\|\langle x\rangle^{3s}\partial^\alpha \left(v^0,   H^0\right)\right\|_0<\infty, \end{align}    i.e.
  \begin{align*}
\|(v-v^{\mm{L}},H-H^{\mm{L}})(t)\|_{k-1} \lesssim\varepsilon{E}^{0}\left(\sqrt{{ {E}}^{0}}
+    \sum_{|\alpha|=k} 
\left\|\langle x\rangle^{3s}\partial^\alpha \left(v^0,   H^0\right)\right\|_0 \right)\mbox{ for any }t\geqslant 0,
\end{align*}}
where $(v^{\mm{L}},H^{\mm{L}}) $ is the unique classical solution to 
 the linear problem 
 \begin{align*}
  &\begin{cases}
   v^\mm{L}_t+\partial_1v^\mm{L}-\varepsilon\nu\Delta v^\mm{L}=0,\\  H^\mm{L}_t+\partial_1H^\mm{L}-\varepsilon\nu\Delta H^\mm{L} =0,\\ 
 \div v^\mm{L}= \div H^\mm{L}=0,\\[1mm]
(v^\mm{L}, H^\mm{L})|_{t=0}=(v^0,H^0).
 \end{cases} 
 \end{align*}
  \item Decay with respect to $\varepsilon$, i.e.
\begin{align*}
|v(x,t)|,\ |H(x,t)|
\lesssim 
  (\langle x+\varepsilon^{-1} \mathbf{e}^1t\rangle^{-s}+\langle x- \varepsilon^{-1} \mathbf{e}^1t
  \rangle^{-s}) \sqrt{{ {E}}^0_{}}
\end{align*}
for any $(x,t)\in \mathbb{R}^n\times \mathbb{R}_0^+$.
In particular,
$$v(x,t),\ H(x,t)\to 0\mbox{ as }\varepsilon\to 0\mbox{ where }t>0.$$ 
If additionally the initial data enjoy the regularity \eqref{25042219501} and have the compact supported in $$B_R(0):=\left\{x\in\mathbb{R}^n\big{|}|x|<R\right\}, $$ the following estimate holds for the case $\nu=0$
 \begin{align*}
|(v,H)(x,t)|\lesssim  \sqrt{\varepsilon}({E}_{}^0)^{\frac{3}{4}} 
\end{align*}
for any $x\in B_R(0)$ and for any $t>2R\varepsilon$.
\end{enumerate}
 \end{thm} 
 
Now we mention Theorem \ref{thsm3}.   Physicists have pointed out that, in the MHD system, a sufficiently
strong magnetic field will reduce the nonlinear interaction and inhibit the formation of strong gradients \cite{kraichnan1965inertial}. In particular, if the following condition of strong field intensity is satisfied
 \begin{align}
 \frac{\mbox{initial perturbation around $(0,\bar{M})$}}{\mbox{field intensity of $\bar{M}$ (denoted by $m$)}}\ll 1, 
 \label{2504061347}
 \end{align}
 where $A\ll B$ means that $A$ is much smaller than $B$ and the impressed magnetic field $\bar{M}$ is a non-zero constant vector, 
 then Bardos--Sulem--Sulem proved that the {ideal} incompressible MHD system (i.e. \eqref{1sa.1} with $\tilde{\mu}=\nu=0$) admits a global solution \cite{MR920153}. The condition \eqref{2504061347} includes the large initial perturbation under the strong field intensity. Later Zhang established the existence result of global solutions for  the 2D incompressible viscous non-resistive MHD system (i.e. \eqref{1sa.1} with $\tilde{\mu}>0$ and $\nu=0$) in $\mathbb{R}^2$ under the condition \eqref{2504061347} by the theory of oscillatory integrals, the method of spectrum analysis and the vanishing mechanism of nonlinear interactions \cite{MR3448784}. Recently, Jiang--Jiang further  extended Zhang's result to the case of 3D periodic domain with the additional asymptotic stability condition of flow maps \cite{MR4641656}, and provided the decay rate of the nonlinear interactions vanishing  with respect to field intensity in Lagrangian coordinates. Recalling the inversely proportional relationship between $\varepsilon$ and $m$, Theorem \ref{thsm3} presents that \emph{the theory of large perturbation solutions for the non-resistive case under strong field intensity can be extended to the case $\tilde{\mu}=\nu>0$}. 
  
 In addition, there are many results of the small Alfv\'en number limit with a 2D limit system for the  \emph{local-in-time} solutions of 3D \emph{non-resistive} MHD equations, see \cite{MR1039472,MR1339828,MR4028271,MR3836806,MR3962511,MR4221327,MR4853427,MR4425021} for examples. However, we further provide the first result of Alfv\'en number limit in Theorem \ref{thsm3} not only  for \emph{global-in-time} solutions, but also for \emph{the resistive case}.

\subsection{An additional result for dissipation vanishing}
It should be noted that the stability estimate in \eqref{2504160945} is also uniform bound with respect to $\nu$, we naturally expect the dissipation vanishing, which reads as follows and the derivation of which will be presented in Section \ref{sec:06} by a standard energy method.
 \begin{thm}\label{2504160948}
	Let $\Lambda^\pm$ be the solutions in Theorem \ref{thm1} with three fixed parameters $n$, $s$,  $k$ and fixed initial data $\Lambda^0_\pm$, which satisfy
\begin{align*} 
 \left\|\left(|\nabla|^{-1}\Lambda^0_+,|\nabla|^{-1}\Lambda^0_-\right)\right\|_0^2+\left\|\langle x\rangle^s\left(\Lambda^0_+, \Lambda^0_-\right)\right\|_0^2\ +\sum_{|\alpha|=1}^{k} 
\left\|\langle x\rangle^{2s}\left(\partial^\alpha\Lambda^0_+, \partial^\alpha\Lambda^0_-\right)\right\|_0^2 \leqslant c_1 \varepsilon^{-\gamma},\end{align*}
 thus  $\Lambda^\pm$ depend on $\nu\geqslant 0$. We denote  $\Lambda^\pm$  by $\Lambda^\pm_{\nu>0}$  resp. $\Lambda^\pm_{\nu=0}$ for the cases $\nu>0$ resp. $\nu=0$, then it holds that
 		\begin{align}
 \label{2504162128}
 		\|\Lambda_{\nu>0}^\pm(t,{{x}})-\Lambda_{\nu=0}^\pm (t,{{x}})\|_{k-1}^2
 	\lesssim \nu  \varepsilon  \mathcal{E}^0_k t
 	e^{ c \varepsilon \sqrt{\mathcal{E}^0}t}\mbox{ for any }t>0.
 		\end{align} 
 \end{thm}
 
 Similarly to Theorem \ref{thsm3}, it is easy to write the result for the non-dissipative limit of the original system \eqref{1sa.1} from the above theorem, and however we omit the details for saving the sink. The rest paper is devoted to the proofs of Theorems \ref{thm1}--\ref{thm3} and \ref{2504160948}.

\section{Existence of solutions with uniform-in-$\varepsilon$ estiamtes}\label{sec2}
As mentioned before, the key proof of Theorem \ref{thm1} is to derive the key \emph{a priori} stability estimate \eqref{2504160945} of solutions $\Lambda^\pm$ for the Cauchy problem \eqref{1.4} defined on $\mathbb{R}^n\times [0,T)$ for any given $T>0$. Next we briefly sketch the derivation of \eqref{2504160945}, and please refer \cite{MR3552009} for the omitted details .

Applying $\partial^\alpha$ to the equations in \eqref{1.4} yields the following systems for ${\Lambda}^\pm$:
\begin{equation}\label{2.1}
\begin{cases}
\partial^\alpha({\Lambda}^\pm_t\mp\partial_1 {\Lambda}^\pm)= \partial^\alpha(\varepsilon\nu \Delta {\Lambda}^\pm-\varepsilon  {\Lambda}^\mp \cdot\nabla {\Lambda}^\pm-\nabla p),\\[1mm] 
\mathrm{div}\partial^\alpha{\Lambda}^\pm=0,
\end{cases}
\end{equation}
where  $0\leqslant |\alpha|\leqslant k$ and $k\geqslant 4$. Similarly, we have
\begin{equation}\label{2saf.1}
\begin{cases}
{|\nabla|^{-1}} ({\Lambda}^\pm_t\mp\partial_1 {\Lambda}^\pm)={ {|\nabla|^{-1}}(\varepsilon\nu\Delta {\Lambda}^\pm}- \varepsilon  {\Lambda}^\mp \cdot\nabla {\Lambda}^\pm-\nabla p),\\[1mm] 
\mathrm{div}{|\nabla|^{-1}} {\Lambda}^\pm=0.
\end{cases}
\end{equation}
Noting that 
\begin{align*}
\partial_t\langle x\pm \mathbf{e}^1t\rangle^{is}\mp\partial_1\langle x\pm \mathbf{e}^1t\rangle^{is} =0,
\end{align*}
where $i=0$, $1$ and 
\begin{align}
\label{2253}
1<2s <4/{3},
\end{align}
thus we can multiply the two equations in $\eqref{2.1}_1$ by the weight functions $\langle x\pm \mathbf{e}^1t\rangle^{ is}$, resp., to obtain
\begin{align}
&\partial_t\left(\langle x\pm \mathbf{e}^1t\rangle^{ is}\partial^\alpha{\Lambda}^\pm\right)\mp\partial_1\left(\langle x\pm\mathbf{e}^1t\rangle^{ is}\partial^\alpha {\Lambda}^\pm\right)\nonumber \\
 & =\langle x\pm \mathbf{e}^1t\rangle^{ is}\partial^\alpha({\varepsilon\nu  \Delta {\Lambda}^\pm}
- \varepsilon( {\Lambda}^{\mp} \cdot\nabla {\Lambda}^\pm)-\nabla p ).\label{2.2} 
\end{align}

Let $\sigma^\pm=\pm x_1-t$ and $q(y)=\int_0^{y}\langle\tau\rangle^{ -2s}{\mm{d}}\tau$. Due to the condition \eqref{2253}, it is easy to check that 
$|q(y)|\lesssim 1$,
\begin{align}
&1\lesssim e^{q(\sigma^\pm)} \leqslant \sup_{(x_1,t)\in \mathbb{R}\times \mathbb{R}_0^+}\{e^{q(\sigma^\pm)}\}=:\tilde{c} \lesssim 1
\label{2504072112}
\end{align}
and  \begin{align}\label{202504171405}
 \partial_te^{q(\sigma^\pm)}\mp \partial_1 e^{q(\sigma^\pm)} =-2e^{q(\sigma^\pm)}\langle x_1\mp t\rangle^{-2s} .
 \end{align}

Taking the $L^2$-inner product of \eqref{2.2}  and $ \langle  {{x}}\pm\mathbf{e}^1 t\rangle^{2{s}} \partial^{ \alpha}  \Lambda_{}^\pm e^{q(\sigma^\pm)}$ for $1\leqslant |\alpha|\leqslant k$ and $i=2$, then summing the identities up for $1\leqslant |\alpha|\leqslant k$, and finally using the integration by parts and the relation \eqref{202504171405}, 
one obtains
 \begin{align}\label{2.3}
 &\frac{1}{2} \frac{\mathrm{d}}{\mathrm{d} t}{\sum_{|\alpha|=1}^k}\int \left|
 \langle x\pm \mathbf{e}^1t\rangle^{2s}{\partial^\alpha}\Lambda^\pm \right|^{2} e^{{q(\sigma^\pm)}}\mm{d}x+{\sum_{|\alpha|=1}^k}\int\frac{\left|\langle x\pm \mathbf{e}^1t\rangle^{2s}{\partial^\alpha}\Lambda^\pm \right|^{2} e^{ {q(\sigma^\pm)} }}{\langle x_1\mp t\rangle^{2s}} \mm{d}x\nonumber \\
 &=I_{\nu,1}^\pm+I_{\mm{N},1}^\pm,
 \end{align}
 where we have defined that 
\begin{align*}
&I_{\nu,1}^\pm:={\varepsilon\nu \sum_{|\alpha|=1}^k  \int\langle x\pm \mathbf{e}^1t\rangle^{4s}\partial^\alpha\Delta {\Lambda}^\pm\cdot  \partial^\alpha\Lambda^\pm e^{q(\sigma^\pm)}{\mm{d}}x},\\
&I_{\mm{N},1}^\pm:= -\sum_{|\alpha|=1}^k\int\langle x\pm \mathbf{e}^1t\rangle^{4s} \partial^\alpha\left(\varepsilon {\Lambda}^\mp \cdot\nabla {\Lambda}^\pm+\nabla p\right)\cdot  \partial^\alpha \Lambda^\pm e^{q(\sigma^\pm)}{\mm{d}}x.
\end{align*}
Similarly, taking the $L^2$-inner product of $\eqref{2.2}$ and $\langle x\pm\mathbf{e}^1t\rangle^{s}\Lambda^\pm e^{q(\sigma^\pm)}$ for  $|\alpha|=0$ and $i=1$, one has 
\begin{align}\label{2.10}
	&\frac{1}{2} \frac{\mathrm{d}}{\mathrm{d} t}\int\left|\langle x\pm \mathbf{e}^1t\rangle^{ s}\Lambda^\pm  \right|^{2}e^{ {q(\sigma^\pm)} }{\mm{d}}x +\int \frac{\left|\langle x\pm \mathbf{e}^1t\rangle^{ s}\Lambda^\pm\right|^2 e^{{q(\sigma^\pm)}}}{\langle x_1\mp t\rangle^{2s}}{\mm{d}}x =I_{\nu,2}^\pm+I_{\mm{N},2}^\pm ,
	\end{align}
	where we have defined that
	\begin{align*}
	&I_{\nu,2}^\pm:=\varepsilon\nu \int\langle x\pm \mathbf{e}^1t\rangle^{2s}\Delta {\Lambda}^\pm\cdot  \Lambda^\pm e^{q(\sigma^\pm)}{\mm{d}}x,\\
&	I_{\mm{N},2}^\pm:=-\int\langle x\pm \mathbf{e}^1t\rangle^{2s}\left(\varepsilon {\Lambda}^\mp \cdot\nabla {\Lambda}^\pm+\nabla p\right)\cdot \Lambda^\pm e^{q(\sigma^\pm)}{\mm{d}}x.
	\end{align*}
In addition, taking the $L^2$-inner product of \eqref{2saf.1}$_1$ and $|\nabla|^{-1}\Lambda^\pm$, and then using the integration by parts and \eqref{2saf.1}$_2$, one gets
	\begin{align}\label{2.13}
	&\frac{1}{2} \frac{\mathrm{d}}{\mathrm{d} t}\int\left||\nabla|^{-1}\Lambda^\pm\right|^{2}\mm{d}x =
I_{\nu,3}^\pm+ I_{\mm{N},3} ^\pm,
	\end{align} 
where we have defined that 
\begin{align*}
&I_{\nu,3}^\pm:=  \varepsilon\nu \int |\nabla|^{-1}\Delta {\Lambda}^\pm\cdot |\nabla|^{-1} \Lambda^\pm  {\mm{d}}x, \\
& I_{\mm{N},3}^\pm:=-\int|\nabla|^{-1}(\varepsilon \Lambda^\mp\cdot  \nabla\Lambda^\pm + \nabla p_{})\cdot|\nabla|^{-1}\Lambda^\pm {\mm{d}}x.
 	\end{align*}      
Multiplying \eqref{2.3},  \eqref{2.10} and \eqref{2.13} by $2(16)^{{k}-|\alpha|}  $,  $ 2(16)^{{k}}   $ and $2\tilde{c}(16)^{k+1}$ (see \eqref{2504072112} for the definition of  $\tilde{c}$), resp., then summing the resulting identities up, and finally integrating them over $(0,t)$, one arrives at 
\begin{align}\label{3.10}
 	& \sup_{0\leqslant \tau\leqslant t}\tilde{\mathcal{E}}^\pm(\tau)\notag + 2\int_0^t  \tilde{\mathcal{W}}^\pm(\tau) \mm{d}\tau\\
 	&	= \tilde{\mathcal{E}}^\pm|_{t=0} +  2\int_0^t\left((16)^{{k}-|\alpha|}(I_{\nu,1}^\pm+I_{\mm{N},1}^\pm)+  (16)^{{k}}   (I_{\nu,2}^\pm+I_{\mm{N},2}^\pm ) +(16)^{k+1}\tilde{c}(I_{\nu,3}^\pm +I_{\mm{N},3}^\pm)\right)\mm{d}\tau,
 	\end{align}
 where we have defined that 
 \begin{align*}
 &\tilde{\mathcal{E}}^\pm(\tau) :=\int  e^{q(\sigma^\pm)}  \Bigg( \sum_{|\alpha|=1}^{{k}}(16)^{{k}-|\alpha|}\big|\langle  {{x}}\pm\mathbf{e}^1  \tau\rangle^{2{s}} \partial^{ \alpha}  \Lambda_{}^\pm\big|^2 \\
   &\qquad \qquad + (16)^{{k}} \big|\langle {{x}}\pm\mathbf{e}^1 \tau\rangle^{{s}} \Lambda_{}^\pm\big|^2+(16)^{k+1}\tilde{c}\mm{sign}(\nu)\left||\nabla|^{-1}\Lambda^\pm\right|^{2}\Bigg)\mm{d}x,\\
 &\tilde{\mathcal{W}}^\pm(\tau):= \int  e^{q(\sigma^\pm)}\Bigg( \sum_{|\alpha|=1}^{{k}}\frac{(16)^{{k}-|\alpha|}\big|\langle  {{x}}\pm\mathbf{e}^1 \tau\rangle^{2{s}} \partial^{\alpha}  \Lambda_{}^\pm\big|^2}{\langle x_1\mp  \tau\rangle^{2{s}}}  +\frac{(16)^{{k}} \big|\langle {{x}}+\mathbf{e}^1 \tau\rangle^{{s}} \Lambda_{}^\pm\big|^2}{\langle x_1\mp\tau\rangle^{2{s}}}  \Bigg)\mm{d}{{x}}.
 	\end{align*}
  
Using the integration by parts, we have
 	\begin{align}\label{3.11}
 	&   (16)^{{k}-|\alpha|}I_{\nu,1}^\pm +  (16)^{{k}}  I_{\nu,2}^\pm +(16)^{k+1}\tilde{c}I_{\nu,3}^\pm 
  = \varepsilon\nu\left( I_{\nu,4}^\pm + I_{\nu,5}^\pm- \tilde{\mathcal{D}} \right) ,
 	\end{align}
 where we have defined that 
 \begin{align*}
 &\tilde{\mathcal{D}}^\pm  :=
 \int\bigg( e^{q(\sigma^\pm)}\Bigg(\sum_{|\alpha|=1}^{{k}}(16)^{{k}-|\alpha|} \Big| \langle{{x}}\pm \mathbf{e}^1\tau\rangle^{2{s}}\nabla \partial^{ \alpha} \Lambda_{}^\pm \big|^2+(16)^{{k}} \big|\langle  {{x}}\pm\mathbf{e}^1  \tau\rangle^{{s}}\nabla \Lambda_{}^\pm \Big|^2\Bigg) \\
 &\qquad \qquad +(16)^{{k}+1} \tilde{c} \left| \Lambda_{}^\pm \right|^2\Bigg)\mm{d}{{x}},\\
 	&I_{\nu,4}^\pm :=-\int \nabla e^{q(\sigma^\pm)}\cdot\Bigg(\sum_{ |\alpha|=1}^{{k}} (16)^{{k}-|\alpha|}\langle{{x}}\pm \mathbf{e}^1\tau\rangle^{4{s}}  \nabla\partial^{\alpha}\Lambda_{}^\pm\cdot \partial^{\alpha} \Lambda_{}^\pm \notag\\
 	& \qquad\qquad+(16)^{{k}}\langle  {{x}}\pm\mathbf{e}^1  \tau\rangle^{2{s}}  \nabla \Lambda_{}^\pm\cdot  \Lambda_{}^\pm \Bigg)\mm{d}{{x}} \notag\\
 	& I_{\nu,5}^\pm:=- \int e^{q(\sigma^\pm)}\Bigg( \sum_{|\alpha|=1}^{{k}}(16)^{{k}-|\alpha|}\nabla\langle{{x}}\pm \mathbf{e}^1\tau\rangle^{4{s}} \cdot \nabla\partial^{ \alpha } \Lambda_{}^\pm\cdot \partial^{\alpha} \Lambda_{}^\pm  \notag\\
 	& \qquad\qquad + (16)^{k}\nabla\langle  {{x}}\pm\mathbf{e}^1  \tau\rangle^{2{s}}\cdot\nabla \Lambda_{}^\pm\cdot \Lambda_{}^\pm \Bigg)\mm{d}{{x}},
 \end{align*}
Exploiting Young's inequality and \eqref{2253}, it is easy to see that 
 	\begin{align}
 	I_{\nu,4}^\pm  	=&\mp  \int \frac{e^{q(\sigma^\pm)}}{\langle x_1\mp \tau\rangle^{2{s}}}\Bigg(\sum_{ |\alpha|=1}^{{k}}(16)^{{k}-|\alpha|}\langle{{x}}\pm \mathbf{e}^1\tau\rangle^{4{s}}  \partial_1\partial^{\alpha } \Lambda_{}^\pm\cdot \partial^{\alpha } \Lambda_{}^\pm \notag\\
 	& + (16)^{{k}}\langle  {{x}}\pm\mathbf{e}^1  \tau\rangle^{2{s}}   \partial_1 \Lambda_{}^\pm\cdot  \Lambda_{}^\pm  \Bigg) \mm{d}{{x}} \notag\\
 	 \leqslant& \int
  e^{q(\sigma^\pm)}\Bigg(\sum_{|\alpha|=1}^{{k}}(16)^{{k}-|\alpha|}\Bigg(\frac{1}{4}\big|\langle{{x}}\pm \mathbf{e}^1\tau\rangle^{2{s}} \partial_1\partial^{\alpha } \Lambda_{}^\pm\big|^2+  \frac{ \big|\langle{{x}}\pm \mathbf{e}^1\tau\rangle^{2{s}}\partial^{\alpha } \Lambda_{}^\pm\big|^2}{\langle x_1\mp\tau\rangle^{2{s}}}\Bigg)\nonumber \\ 
 	 & +(16)^{{k}}\Bigg(\frac{1}{4} \big|\langle  {{x}}\pm\mathbf{e}^1  \tau\rangle^{{s}} \partial_1 \Lambda_{}^\pm \big|^2+\frac{ \big|\langle  {{x}}\pm\mathbf{e}^1  \tau\rangle^{{s}}  \Lambda_{}^\pm\big|^2}{\langle x_1\mp \tau\rangle^{2{s}}}\Bigg)\Bigg)\mm{d}{{x}} \label{250417}
 	\end{align}
 and
 	\begin{align}
 	 I_{\nu,5}^\pm   =&-{s} \int e^{q(\sigma^\pm)} \Bigg(4\sum_{|\alpha|=1}^{{k}}(16)^{{k}-|\alpha|}\langle{{x}}\pm \mathbf{e}^1\tau\rangle^{4{s}-2}({{x}}\pm \mathbf{e}^1\tau) \cdot \nabla\partial^{\alpha } \Lambda_{}^\pm\cdot\partial^{\alpha } \Lambda_{}^\pm \notag\\
 	&  + 2(16)^{{k}} \langle  {{x}}\pm\mathbf{e}^1  \tau\rangle^{2{s}-2}({{x}}\pm\mathbf{e}^1  \tau)\cdot\nabla \Lambda_{}^\pm\cdot  \Lambda_{}^\pm \Bigg)\mm{d}{{x}} \notag\\
 	\leqslant& (16)^{{k}}\int e^{q(\sigma^\pm)} \Bigg(\sum_{|\alpha|=1}^{{k}}(16)^{ -|\alpha|}\Bigg(\frac{1}{4}\big| \langle{{x}}\pm \mathbf{e}^1\tau\rangle^{2{s}}\nabla\partial^{ \alpha} \Lambda_{}^\pm \big|^2+\frac{15}{2}  \big|\langle{{x}}\pm \mathbf{e}^1\tau\rangle^{2{s}-1}\partial^{  \alpha} \Lambda_{}^\pm\big|^2\Bigg) \notag\\ &  +\frac{1}{4}\big|\langle  {{x}}\pm\mathbf{e}^1  \tau\rangle^{{s}} \nabla \Lambda_{}^\pm \big|^2 + \frac{15}{2}\big| \Lambda_{}^\pm\big|^2\Bigg) \mm{d}{{x}} .\label{2504171}
 	\end{align}
 
Making use of \eqref{2504072112}  and \eqref{3.11}--\eqref{2504171}, we infer from \eqref{3.10} that
 \begin{align}\label{3.10saf}
 \sup_{0\leqslant \tau< t}\tilde{\mathcal{E}}^\pm(\tau) +\int_0^t\left((1-\varepsilon\nu)  \tilde{\mathcal{W}}^\pm(\tau) +\varepsilon\nu \tilde{\mathcal{D}}^\pm(\tau)\right)\mm{d}\tau \lesssim  \tilde{\mathcal{E}}^\pm|_{t=0} +  \int_0^t\sum_{i=1}^3\left|I_{\mm{N},i}^\pm (\tau)\right|\mm{d}\tau,
 	\end{align}
 where $t\in [0,T]$.
  In addition, the first author and Lei in \cite{MR3552009} have proved that
 \begin{align}
  \int_0^t\sum_{i=1}^3\left|I_{\mm{N},i}^\pm (\tau)\right| \mm{d}\tau\lesssim  \sup_{0\leqslant \tau<t}\sqrt{\mathcal{E}_k(\tau)  } \int_0^t \mathcal{W}_k (\tau) \mm{d}\tau\mbox{ for the case }\varepsilon=1.
  \label{2504221914}
  \end{align}
 Recalling the expression of $p$, i.e.
 	\begin{align*} 
 	p= \varepsilon  (-\Delta)^{-1}\sum_{1\leqslant i, j\leqslant n}\partial_{x_i}\partial_{x_j} \left( \Lambda^-_i  \Lambda^+_j\right).
 	\end{align*}  
 thus, it is easy to see from the derivation of \eqref{2504221914} that
   \begin{align*}
  \int_0^t\sum_{i=1}^3\left|I_{\mm{N},i}^\pm (\tau)\right| \mm{d}\tau\lesssim \varepsilon  \sup_{0\leqslant \tau<t}\sqrt{\mathcal{E}_k(\tau)  } \int_0^t \mathcal{W}_k (\tau) \mm{d}\tau\mbox{ for any }\varepsilon>0.
  \end{align*}
 Putting the above estimate into \eqref{3.10saf}, and then using \eqref{2504072112}, we immediately obtain
   \begin{align}
 \sup_{0\leqslant \tau< t} {\mathcal{E}}_k(\tau) +\int_0^t\left((1- \varepsilon\nu)  {\mathcal{W}}_k(\tau) +\varepsilon\nu  {\mathcal{D}}_k(\tau)\right)\mm{d}\tau \lesssim{\mathcal{E}}^0     +  \varepsilon  \sup_{0\leqslant \tau<t}\sqrt{\mathcal{E}_k(\tau)  } \int_0^t \mathcal{W}_k (\tau) \mm{d}\tau 
 \label{2504082303}
 	\end{align}
 for any $t\in [0,T]$.
 Hence solutions $\Lambda^\pm$ satisfy \eqref{2504161849}.
Consequently, under the \emph{a priori}  assumption \eqref{2504171021}, we further have \eqref{2504161849}, which, together with a local well-posedness result on the Cauchy problem \eqref{1.4}, immediately yields Theorem \ref{thm1}.	 

\section{An error estimate for nonlinear interactions vanishing}\label{sec3}

This section is dedicated to the proof of Theorem \ref{thm2}. 
Obviously, the linear problem \eqref{thm2} admits unique solutions $\Lambda^\pm_{\mm{L}}$, which enjoy the estimate
\begin{align}
&	\mathcal{E}^\mm{L}_k(t) +\int_0^t\left(\mathcal{W}^\mm{L}_k(\tau)+\varepsilon\nu \mathcal{D}^\mm{L}_k(\tau) \right){\mm{d}}\tau \lesssim \mathcal{E}^0_k  
  \label{3.1}
\end{align}
and, by the regularity condition \eqref{2504221950},
 \begin{align}
&	 \sum_{|\alpha|=  k}\left\|\left(\langle x+ \mathbf{e}^1t\rangle^{3s}\partial^\alpha\Lambda^+_\mm{L},\langle x- \mathbf{e}^1t\rangle^{3s}\partial^\alpha\Lambda^-_\mm{L} \right)\right\|_0^2\nonumber \\ 
&+\int_0^t\sum_{|\alpha|= k}\left\|\left(\frac{\langle x+ \mathbf{e}^1t\rangle^{3s}\partial^\alpha\Lambda^+_\mm{L}}{\langle x_1- t\rangle^{s}},\frac{\langle x- \mathbf{e}^1t\rangle^{3s}\partial^\alpha\Lambda^-_\mm{L}}{\langle x_1+ t\rangle^{s}}\right)\right\|_0^2 {\mm{d}}\tau\nonumber \\
&\lesssim \sum_{|\alpha|=k} 
\left\|\langle x\rangle^{3s}\left(\partial^\alpha\Lambda^0_+, \partial^\alpha\Lambda^0_-\right)\right\|_0^2 
  \label{3.11}
\end{align}
where $\mathcal{E}^\mm{L}_k(t)$ resp. $\mathcal{D}^\mm{L}_k(t)$ are defined by $\mathcal{E}_k(t)$ resp. $\mathcal{D}_k(t)$ with $\Lambda^\pm_\mm{L}$ in place of $\Lambda^\pm$.

Let $\Lambda^\pm_{\mm{N},\mm{L}}=\Lambda^\pm-\Lambda^\pm_\mm{L}$, then the error functions $\Lambda^\pm_{\mm{N},\mm{L}}$ satisfy
\begin{equation}\label{3.2}
\begin{cases}
\partial_t{\Lambda}^{\pm}_{\mm{N},\mm{L}}\mp\partial_1 {\Lambda}^\pm_{\mm{N},\mm{L}}-\varepsilon\nu\Delta {\Lambda}_{\mm{N},\mm{L}}^\pm=-\varepsilon {\Lambda}^\mp \cdot\nabla {\Lambda}^\pm-\nabla p ,\\[1mm] 
\mathrm{div}{\Lambda}^\pm_{\mm{N},\mm{L}}=0,\\[1mm]
\Lambda^\pm_{\mm{N},\mm{L}}|_{t=0}=0.
\end{cases}
\end{equation}
In adaption, making use of the imbedding inequality with weight in \cite[Lemma 2.1]{MR3780143},  the divergence-free conditions of $\nabla {\Lambda}^\pm$ in \eqref{1.4}$_2$, the integration by parts and Young's inequality, we have
\begin{align}
&\int {\Lambda}^\mp \cdot\nabla\partial^\alpha {\Lambda}^\pm\cdot  \partial^\alpha \Lambda^\pm_{\mm{N},\mm{L}}\mm{d}x\nonumber \\
&= \int {\Lambda}^\mp \cdot\nabla  \partial^\alpha {\Lambda}^\pm_{\mm{N},\mm{L}} \cdot  \partial^\alpha \Lambda^\pm_{\mm{N},\mm{L}}\mm{d}x+
 \int {\Lambda}^\mp \cdot\nabla \partial^\alpha{\Lambda}^\pm_{\mm{L}}\cdot  \partial^\alpha \Lambda^\pm_{\mm{N},\mm{L}}\mm{d}x\nonumber \\
 & = 
 \int {\Lambda}^\mp \cdot\nabla \partial^\alpha{\Lambda}^\pm_{\mm{L}}\cdot  \partial^\alpha \Lambda^\pm_{\mm{N},\mm{L}}\mm{d}x  \lesssim    \left\|\langle x\pm \mathbf{e}^1t\rangle^{3s}\nabla\partial^\alpha\Lambda^{\pm}_{\mm{L}}\right\|_0 ( \mathcal{W}_k  +\mathcal{W}^\mm{L}_k ), \nonumber 
\end{align}
where $|\alpha|=k-1$. 

Keeping the above estimate in mind, and then following the same process as in the derivation of \eqref{2504082303},  we can deduce from \eqref{3.2} that
\begin{align}
  &   \mathcal{E}_{k-1}^{\mm{N},\mm{L}}(t) +\int_0^t({\mathcal{W}}_{k-1}^{\mm{N},\mm{L}}(\tau)+\varepsilon\nu  
{\mathcal{D}}_{k-1}^{\mm{N},\mm{L}}(\tau)){\mm{d}}\tau \nonumber \\
  &\lesssim \varepsilon  \sup_{0\leqslant \tau< t} \bigg( \mathcal{E}_{k-1}^{\mm{N},\mm{L}}(\tau)+\mathcal{E}_{k-1}^{\mm{L}}(\tau) \nonumber 
 \\
  &\qquad +\left\|\left(\langle x+ \mathbf{e}^1t\rangle^{3s}\nabla\partial^\alpha\Lambda^+_\mm{L},\langle x- \mathbf{e}^1t\rangle^{3s}\nabla\partial^\alpha\Lambda^-_\mm{L} \right)\right\|_0^2 \bigg)^{1/2} \int_0^t 
   ( \mathcal{W}_{k-1}^{\mm{N},\mm{L}}(\tau) +\mathcal{W}^\mm{L}_{k-1}(\tau))\mm{d}\tau \nonumber
 \end{align}  for any $t>0$.
Finally, inserting \eqref{1.12}, \eqref{3.1} and \eqref{3.11} into the above estimate   yields the desired estimate \eqref{1.12}. This completes the proof of Theorem \ref{thm2}.
	
\section{Decay of solutions}\label{sec4}
	
Now we turn to the verification of Theorem \ref{thm3}. 
Exploiting the Sobolev embedding inequality $H^2\hookrightarrow L^\infty$, we have 
\begin{align}
|\langle x\pm \mathbf{e}^1t\rangle^s\Lambda^\pm|\lesssim  \sum_{|\alpha|=0}^2\| \partial^\alpha(\langle x\pm \mathbf{e}^1t\rangle^s \Lambda^\pm)\|_0\lesssim 
 \sum_{|\alpha|=0}^2\| \langle x\pm \mathbf{e}^1t\rangle^s{\partial^\alpha}\Lambda^\pm\|_0\lesssim \sqrt{\mathcal{E}^{0}_{k}},
 \label{2504140913}
\end{align}
which yields
\begin{align*}
|\Lambda^\pm(x,t)|\lesssim    \langle x\pm \mathbf{e}^1t\rangle^{-s}\sqrt{\mathcal{E}^{0}_{k}}\mbox{ for any }(x,t)\in \mathbb{R}^n\times \mathbb{R}_0^+. 
\end{align*}
Hence \eqref{2504140938} holds for $\nu\geqslant 0$. Next we shall derive {\eqref{2504140939}} for the case $\nu=0$.

The solutions of the linear Cauchy problem \eqref{1.11} are given by the formulas \cite{MR1625845}:
\begin{align}
\label{2022501021627}
\Lambda^\pm_\mm{L}=\Lambda_\pm^0(x\pm \mathbf{e}^1t).
\end{align}
Similarly to \eqref{2504140913}, it holds that 
$$
|\langle x \rangle^s \Lambda^0_\pm|\lesssim  \sum_{|\alpha|=0}^2\| \partial^\alpha(\langle x  \rangle^s \Lambda^0_\pm )\|_0\lesssim 
 \sum_{|\alpha|=0}^2\| \langle x \rangle^s \partial^\alpha\Lambda^0_\pm\|_0\lesssim  \sqrt{\mathcal{E}^{0}_k},
$$
which yields
\begin{align}
\label{2501021627}\Lambda^0_\pm\lesssim  \sqrt{\mathcal{E}^{0}_k} \langle x \rangle^{-s}\mbox{ for }x\neq 0. 
\end{align}
It follows from \eqref{2022501021627} and \eqref{2501021627} that
\begin{align*} 
\Lambda^\pm_\mm{L}(x,t)\lesssim \sqrt{\mathcal{E}^{0}_{k}} \langle  x\pm \mathbf{e}^1t  \rangle^{-s}\mbox{ for any }(x,t)\in \mathbb{R}^n\times \mathbb{R}_0^+.
\end{align*}

If additionally ${\rm supp}\Lambda_\pm^0\subset B_R(0)$, it is easy to see from \eqref{2022501021627} that 
\begin{align*}
\Lambda^\pm_\mm{L}\equiv0\ {\rm for}\ |x\pm \mathbf{e}^1t|>R,
\end{align*}
which implies that, for any $x\in B_R(0)$ and for any $t>2R$, 
\begin{align}\label{4.1}
\Lambda^\pm_\mm{L}\equiv0.
\end{align}
Exploiting \eqref{1.12} and \eqref{4.1}, we conclude that 
\begin{align*}
|\Lambda^\pm(x,t)|\leqslant|\Lambda^\pm_\mm{L}(x,t)|+|\Lambda^\pm_{\mm{N},\mm{L}}(x,t)|
\lesssim  \varepsilon^{\frac{1}{2}}(\mathcal{E}_{k}^0)^{\frac{3}{4}} 
\end{align*}for any $x\in B_R(0)$ and for any $t>2R$. Therefore {\eqref{2504140939}} holds for the case $\nu=0$. 
This completes the proof of Theorem \ref{thm3}.

 \section{Non-dissipative limit}\label{sec:06}
This section is devoted to the proof of non-dissipative limit stated in Theorem \ref{2504160948}.
To this purpose, we define that $\Lambda_{\nu,0}^\pm:=\Lambda_{\nu>0}^\pm-\Lambda_{\nu=0}^\pm$ and  $p_{\nu,0}:=p_{\nu>0}-p_{\nu=0}$, where $p_{\nu>0}$ resp.  $p_{\nu=0}$ are the pressure functions for the cases $\nu>0$ resp. $\nu=0$. Thus we obtain the following systems for $\Lambda_{\nu,0}^\pm$.
\begin{align}\label{3.15}
\begin{cases}
 	\partial_t \Lambda_{\nu,0}^\pm\mp\partial_1\Lambda_{\nu,0}^\pm+\nabla p_{\nu,0}-\varepsilon\nu\Delta \Lambda_{\nu,0}^\pm-\varepsilon\nu\Delta \Lambda_{\nu=0}^\pm\\
 	=-\varepsilon\Big(\Lambda_{\nu,0}^{\mp}\cdot\nabla \Lambda_{\nu,0}^\pm+\Lambda_{\nu,0}^{\mp}\cdot\nabla  \Lambda_{\nu=0}^\pm+\Lambda_{\nu=0}^{\mp}\cdot\nabla \Lambda_{\nu,0}^\pm\Big),\\ 
 	\mm{div} \Lambda_{\nu,0}^\pm=0.
 	\end{cases}
 	\end{align}
 
 Let $\alpha$ satisfy $|\alpha|\leqslant {k}-1$.
 	Applying $\partial^\alpha$ to the above systems, then taking the inner product of the resulting 
 identities and $2\partial^\alpha \Lambda_{\nu,0}^\pm$ in $L^2$, resp. and finally using the embedding inequality of $H^2\hookrightarrow L^\infty $,
  the divergence-free conditions in \eqref{3.15}$_2$, the integration by parts, one  gets
 	\begin{align*}
 	&\frac{\mm{d}}{\mm{d}t}\| \partial^\alpha \Lambda_{\nu,0}^\pm\|^2_0+2\varepsilon\nu\|\nabla \partial^\alpha \Lambda_{\nu,0}^\pm\|^2_0\notag\\
 	&=-2\varepsilon\nu\int \Delta\partial^\alpha \Lambda_{\nu=0}^\pm\cdot \partial^\alpha \Lambda_{\nu,0}^\pm\mm{d}{{x}}-\varepsilon\int \partial^\alpha\Big(\Lambda_{\nu,0}^\mp\cdot\nabla \Lambda_{\nu,0}^\pm+\Lambda_{\nu,0}^\mp\cdot\nabla \Lambda_{\nu=0}^\pm\notag\\
 	&\quad +\Lambda_{\nu=0}^\mp\cdot\nabla \Lambda_{\nu,0}^\pm\Big)\cdot\partial^\alpha \Lambda_{\nu,0}^\pm\mm{d}{{x}}\notag\\
 	&\leqslant 2\varepsilon\nu \|\nabla \partial^\alpha \Lambda_{\nu,0}^\pm\|_0 \|\nabla \partial^\alpha \Lambda_{\nu=0}^\pm\|_0 +c\varepsilon\Big(\|(\Lambda_{\nu,0}^+,\Lambda_{\nu,0}^-)\|_{[{k}/2]}+\|(\Lambda_{\nu=0}^+,
 \Lambda_{\nu=0}^- )\|_{k}\Big) \|\Lambda_{\nu,0}^\pm\|_{k-1}^2,
 	\end{align*}
which, together with \eqref{2504160945} and Young inequality, implies
 	\begin{align*}
 	\frac{\mm{d}}{\mm{d}t} \|  \Lambda_{\nu,0}^\pm\|
 	_{k-1}^2
 	&\lesssim \varepsilon\Big(\|(\Lambda_{\nu,0}^+,\Lambda_{\nu,0}^-)\|_{[{k}/2]}+\|(\Lambda_{\nu=0}^+,\Lambda_{\nu=0}^- )\|_{k}\Big) \|  \Lambda_{\nu,0}^\pm\|
 	_{k-1}^2+\nu\varepsilon\|   \Lambda_{\nu=0}^\pm\|_k^2 \\
 	& \lesssim   \varepsilon \sqrt{ \mathcal{E}^0_k}\|  \Lambda_{\nu,0}^\pm\|
 	_{k-1}^2+\nu \varepsilon  \mathcal{E}^0_k .
 	\end{align*}
Applying Gronwall's inequality to the above inequality, we arrive at
 	\begin{align*}
 	\|  \Lambda_{\nu,0}^\pm(t)\|^2
 	_{k-1} \lesssim \int_0^t  \nu \varepsilon  \mathcal{E}^0_k \mm{d}\tau e^{c\int_0^t   \varepsilon \sqrt{\mathcal{E}^0_k}\mm{d}\tau}
 	\lesssim \nu  \varepsilon  \mathcal{E}^0_k t
 	e^{ c \varepsilon \sqrt{\mathcal{E}^0_k}t},
 	\end{align*} 
 which yields \eqref{2504162128}. This completes the proof of Theorem \ref{2504160948}.

\vspace{4mm}
\noindent\textbf{Acknowledgements.} The research of Yuan Cai was supported by NSFC (Grant No. 12201122),     the research of Xiufang Cui by NSFC (Grant No. 12401278), and the research of Fei Jiang by NSFC (Grant Nos. 12371233 and 12231016), and the NSF of Fujian Province of China (Grant  Nos. 2024J011011 and 2022J01105), Fujian Alliance of Mathematics (Grant No. 2025SXLMMS01) and the Central Guidance on Local Science and Technology Development Fund of Fujian Province (Grant No. 2023L3003). The third author thank to Dr. Jiawei Wang, Prof. Jin Xu, and Prof. Yi Zhou for discussion on the small Alfv\'en number limit.


	\renewcommand\refname{References}
	\renewenvironment{thebibliography}[1]{%
		\section*{\refname}
		\list{{\arabic{enumi}}}{\def\makelabel##1{\hss{##1}}\topsep=0mm
			\parsep=0mm
			\partopsep=0mm\itemsep=0mm
			\labelsep=1ex\itemindent=0mm
			\settowidth\labelwidth{\small[#1]}
			\leftmargin\labelwidth \advance\leftmargin\labelsep
			\advance\leftmargin -\itemindent
			\usecounter{enumi}}\small
		\def\newblock{\ }
		\sloppy\clubpenalty4000\widowpenalty4000
		\sfcode`\.=1000\relax}{\endlist}
	\bibliographystyle{model1b-num-names}
	\bibliography{refs}

\begin{thebibliography}{30}
\expandafter\ifx\csname natexlab\endcsname\relax\def\natexlab#1{#1}\fi
\providecommand{\bibinfo}[2]{#2}
\ifx\xfnm\relax \def\xfnm[#1]{\unskip,\space#1}\fi
\bibitem[{Abidi and Zhang(2017)}]{MR3666563}
\bibinfo{author}{H.~Abidi}, \bibinfo{author}{P.~Zhang}, \bibinfo{title}{On the
  global solution of a 3-{D} {MHD} system with initial data near equilibrium},
  \bibinfo{journal}{Comm. Pure Appl. Math.} \bibinfo{volume}{70}
  (\bibinfo{year}{2017}) \bibinfo{pages}{1509--1561}.
\bibitem[{Alinhac(2001)}]{MR1856402}
\bibinfo{author}{S.~Alinhac}, \bibinfo{title}{The null condition for
  quasilinear wave equations in two space dimensions {I}},
  \bibinfo{journal}{Invent. Math.} \bibinfo{volume}{145} (\bibinfo{year}{2001})
  \bibinfo{pages}{597--618}.
\bibitem[{Bardos et~al.(1988)Bardos, Sulem and Sulem}]{MR920153}
\bibinfo{author}{C.~Bardos}, \bibinfo{author}{C.~Sulem}, \bibinfo{author}{P.L.
  Sulem}, \bibinfo{title}{Longtime dynamics of a conductive fluid in the
  presence of a strong magnetic field}, \bibinfo{journal}{Trans. Amer. Math.
  Soc.} \bibinfo{volume}{305} (\bibinfo{year}{1988}) \bibinfo{pages}{175--191}.
\bibitem[{Browning and Kreiss(1982)}]{MR665380}
\bibinfo{author}{G.~Browning}, \bibinfo{author}{H.O. Kreiss},
  \bibinfo{title}{Problems with different time scales for nonlinear partial
  differential equations}, \bibinfo{journal}{SIAM J. Appl. Math.}
  \bibinfo{volume}{42} (\bibinfo{year}{1982}) \bibinfo{pages}{704--718}.
\bibitem[{Cai and Lei(2018)}]{MR3780143}
\bibinfo{author}{Y.~Cai}, \bibinfo{author}{Z.~Lei}, \bibinfo{title}{Global
  well-posedness of the incompressible magnetohydrodynamics},
  \bibinfo{journal}{Arch. Ration. Mech. Anal.} \bibinfo{volume}{228}
  (\bibinfo{year}{2018}) \bibinfo{pages}{969--993}.
\bibitem[{Cheng et~al.(2018)Cheng, Ju and Schochet}]{MR3803773}
\bibinfo{author}{B.~Cheng}, \bibinfo{author}{Q.~Ju},
  \bibinfo{author}{S.~Schochet}, \bibinfo{title}{Three-scale singular limits of
  evolutionary {PDE}s}, \bibinfo{journal}{Arch. Ration. Mech. Anal.}
  \bibinfo{volume}{229} (\bibinfo{year}{2018}) \bibinfo{pages}{601--625}.
\bibitem[{Cheng et~al.(2021)Cheng, Ju and Schochet}]{MR4221327}
\bibinfo{author}{B.~Cheng}, \bibinfo{author}{Q.~Ju},
  \bibinfo{author}{S.~Schochet}, \bibinfo{title}{Convergence rate estimates for
  the low {M}ach and {A}lfv\'en number three-scale singular limit of
  compressible ideal magnetohydrodynamics}, \bibinfo{journal}{ESAIM Math.
  Model. Numer. Anal.} \bibinfo{volume}{55} (\bibinfo{year}{2021})
  \bibinfo{pages}{S733--S759}.
\bibitem[{Evans(1998)}]{MR1625845}
\bibinfo{author}{L.C. Evans}, \bibinfo{title}{Partial differential equations},
  volume~\bibinfo{volume}{19} of \textit{\bibinfo{series}{Graduate Studies in
  Mathematics}}, \bibinfo{publisher}{American Mathematical Society, Providence,
  RI}, \bibinfo{year}{1998}.
\bibitem[{Fefferman et~al.(2017)Fefferman, McCormick, Robinson and
  Rodrigo}]{MR3590662}
\bibinfo{author}{C.L. Fefferman}, \bibinfo{author}{D.S. McCormick},
  \bibinfo{author}{J.C. Robinson}, \bibinfo{author}{J.L. Rodrigo},
  \bibinfo{title}{Local existence for the non-resistive {MHD} equations in
  nearly optimal {S}obolev spaces}, \bibinfo{journal}{Arch. Ration. Mech.
  Anal.} \bibinfo{volume}{223} (\bibinfo{year}{2017})
  \bibinfo{pages}{677--691}.
\bibitem[{Goto(1990)}]{MR1039472}
\bibinfo{author}{S.~Goto}, \bibinfo{title}{Singular limit of the incompressible
  ideal magneto-fluid motion with respect to the {A}lfv\'en number},
  \bibinfo{journal}{Hokkaido Math. J.} \bibinfo{volume}{19}
  (\bibinfo{year}{1990}) \bibinfo{pages}{175--187}.
\bibitem[{He et~al.(2018)He, Xu and Yu}]{MR3738384}
\bibinfo{author}{L.B. He}, \bibinfo{author}{L.~Xu}, \bibinfo{author}{P.~Yu},
  \bibinfo{title}{On global dynamics of three dimensional magnetohydrodynamics:
  nonlinear stability of {A}lfv\'en waves}, \bibinfo{journal}{Ann. PDE}
  \bibinfo{volume}{4} (\bibinfo{year}{2018}) \bibinfo{pages}{Paper No. 5, 105}.
\bibitem[{Jiang and Jiang(2023)}]{MR4641656}
\bibinfo{author}{F.~Jiang}, \bibinfo{author}{S.~Jiang}, \bibinfo{title}{On
  magnetic inhibition theory in 3{D} non-resistive magnetohydrodynamic fluids:
  global existence of large solutions}, \bibinfo{journal}{Arch. Ration. Mech.
  Anal.} \bibinfo{volume}{247} (\bibinfo{year}{2023}) \bibinfo{pages}{Paper No.
  96, 35}.
\bibitem[{Jiang et~al.(2019)Jiang, Ju and Xu}]{MR4028271}
\bibinfo{author}{S.~Jiang}, \bibinfo{author}{Q.~Ju}, \bibinfo{author}{X.~Xu},
  \bibinfo{title}{Small {A}lfv\'{e}n number limit for incompressible
  magneto-hydrodynamics in a domain with boundaries}, \bibinfo{journal}{Sci.
  China Math.} \bibinfo{volume}{62} (\bibinfo{year}{2019})
  \bibinfo{pages}{2229--2248}.
\bibitem[{Ju et~al.(2019)Ju, Schochet and Xu}]{MR3962511}
\bibinfo{author}{Q.~Ju}, \bibinfo{author}{S.~Schochet},
  \bibinfo{author}{X.~Xu}, \bibinfo{title}{Singular limits of the equations of
  compressible ideal magneto-hydrodynamics in a domain with boundaries},
  \bibinfo{journal}{Asymptot. Anal.} \bibinfo{volume}{113}
  (\bibinfo{year}{2019}) \bibinfo{pages}{137--165}.
\bibitem[{Ju et~al.(2024)Ju, Wang and Xu}]{MR4645635}
\bibinfo{author}{Q.~Ju}, \bibinfo{author}{J.~Wang}, \bibinfo{author}{X.~Xu},
  \bibinfo{title}{Small {A}lfv\'{e}n number limit for shallow water
  magnetohydrodynamics}, \bibinfo{journal}{J. Math. Anal. Appl.}
  \bibinfo{volume}{531} (\bibinfo{year}{2024}) \bibinfo{pages}{Paper No.
  127773, 23}.
\bibitem[{Ju and Xu(2018)}]{MR3836806}
\bibinfo{author}{Q.~Ju}, \bibinfo{author}{X.~Xu}, \bibinfo{title}{Small
  {A}lfv\'en number limit of the plane magnetohydrodynamic flows},
  \bibinfo{journal}{Appl. Math. Lett.} \bibinfo{volume}{86}
  (\bibinfo{year}{2018}) \bibinfo{pages}{77--82}.
\bibitem[{Klainerman and Majda(1981)}]{MR615627}
\bibinfo{author}{S.~Klainerman}, \bibinfo{author}{A.~Majda},
  \bibinfo{title}{Singular limits of quasilinear hyperbolic systems with large
  parameters and the incompressible limit of compressible fluids},
  \bibinfo{journal}{Comm. Pure Appl. Math.} \bibinfo{volume}{34}
  (\bibinfo{year}{1981}) \bibinfo{pages}{481--524}.
\bibitem[{Kraichnan(1965)}]{kraichnan1965inertial}
\bibinfo{author}{R.H. Kraichnan}, \bibinfo{title}{Inertial-range spectrum of
  hydromagnetic turbulence}, \bibinfo{journal}{Physics of Fluids}
  \bibinfo{volume}{8} (\bibinfo{year}{1965}) \bibinfo{pages}{1385}.
\bibitem[{Kuku\v{c}ka(2011)}]{MR2805860}
\bibinfo{author}{P.~Kuku\v{c}ka}, \bibinfo{title}{Singular limits of the
  equations of magnetohydrodynamics}, \bibinfo{journal}{J. Math. Fluid Mech.}
  \bibinfo{volume}{13} (\bibinfo{year}{2011}) \bibinfo{pages}{173--189}.
\bibitem[{Lei(2016)}]{MR3552009}
\bibinfo{author}{Z.~Lei}, \bibinfo{title}{Global well-posedness of
  incompressible elastodynamics in two dimensions}, \bibinfo{journal}{Comm.
  Pure Appl. Math.} \bibinfo{volume}{69} (\bibinfo{year}{2016})
  \bibinfo{pages}{2072--2106}.
\bibitem[{Majda(1984)}]{MR748308}
\bibinfo{author}{A.~Majda}, \bibinfo{title}{Compressible fluid flow and systems
  of conservation laws in several space variables}, volume~\bibinfo{volume}{53}
  of \textit{\bibinfo{series}{Applied Mathematical Sciences}},
  \bibinfo{publisher}{Springer-Verlag, New York}, \bibinfo{year}{1984}.
\bibitem[{Rubino(1995)}]{MR1339828}
\bibinfo{author}{B.~Rubino}, \bibinfo{title}{Singular limits in the data space
  for the equations of magneto-fluid dynamics}, \bibinfo{journal}{Hokkaido
  Math. J.} \bibinfo{volume}{24} (\bibinfo{year}{1995})
  \bibinfo{pages}{357--386}.
\bibitem[{Schochet(1986)}]{MR871107}
\bibinfo{author}{S.~Schochet}, \bibinfo{title}{Symmetric hyperbolic systems
  with a large parameter}, \bibinfo{journal}{Comm. Partial Differential
  Equations} \bibinfo{volume}{11} (\bibinfo{year}{1986})
  \bibinfo{pages}{1627--1651}.
\bibitem[{Schochet(1988)}]{MR957005}
\bibinfo{author}{S.~Schochet}, \bibinfo{title}{Asymptotics for symmetric
  hyperbolic systems with a large parameter}, \bibinfo{journal}{J. Differential
  Equations} \bibinfo{volume}{75} (\bibinfo{year}{1988})
  \bibinfo{pages}{1--27}.
\bibitem[{Schochet(1994)}]{MR1303036}
\bibinfo{author}{S.~Schochet}, \bibinfo{title}{Fast singular limits of
  hyperbolic {PDE}s}, \bibinfo{journal}{J. Differential Equations}
  \bibinfo{volume}{114} (\bibinfo{year}{1994}) \bibinfo{pages}{476--512}.
\bibitem[{Sermange and Temam(1983)}]{MR716200}
\bibinfo{author}{M.~Sermange}, \bibinfo{author}{R.~Temam}, \bibinfo{title}{Some
  mathematical questions related to the {MHD} equations},
  \bibinfo{journal}{Comm. Pure Appl. Math.} \bibinfo{volume}{36}
  (\bibinfo{year}{1983}) \bibinfo{pages}{635--664}.
\bibitem[{Wang and Xu(2025)}]{MR4853427}
\bibinfo{author}{X.~Wang}, \bibinfo{author}{X.~Xu}, \bibinfo{title}{Singular
  limits of compressible viscous {MHD} system with vertical magnetic field},
  \bibinfo{journal}{J. Differential Equations} \bibinfo{volume}{425}
  (\bibinfo{year}{2025}) \bibinfo{pages}{470--505}.
\bibitem[{Wei and Zhang(2017)}]{MR3678491}
\bibinfo{author}{D.~Wei}, \bibinfo{author}{Z.~Zhang}, \bibinfo{title}{Global
  well-posedness of the {MHD} equations in a homogeneous magnetic field},
  \bibinfo{journal}{Anal. PDE} \bibinfo{volume}{10} (\bibinfo{year}{2017})
  \bibinfo{pages}{1361--1406}.
\bibitem[{Zhang(2022)}]{MR4425021}
\bibinfo{author}{S.~Zhang}, \bibinfo{title}{Singular limit of the nonisentropic
  compressible ideal {MHD} equations in a domain with boundary},
  \bibinfo{journal}{Appl. Anal.} \bibinfo{volume}{101} (\bibinfo{year}{2022})
  \bibinfo{pages}{2596--2615}.
\bibitem[{Zhang(2016)}]{MR3448784}
\bibinfo{author}{T.~Zhang}, \bibinfo{title}{Global solutions to the 2{D}
  viscous, non-resistive {MHD} system with large background magnetic field},
  \bibinfo{journal}{J. Differential Equations} \bibinfo{volume}{260}
  (\bibinfo{year}{2016}) \bibinfo{pages}{5450--5480}.

\end{thebibliography}
\end{CJK*}
\end{document}